\RequirePackage{fix-cm}
\documentclass[smallextended]{svjour3}       
\smartqed  
\usepackage{graphicx}

\usepackage{amsmath}
\usepackage{amssymb}
\usepackage{epsfig}
\usepackage{colordvi}
\usepackage{graphics}
\usepackage{color}

\usepackage{cite}
\usepackage{rotating}
\usepackage{booktabs}
\usepackage{lmodern}

 \usepackage{mathptmx}      
 \journalname{Journal of Elasticity}

\newcommand{\bfb}{{\bf b}}%
\newcommand{\bfn}{{\bf n}}%
\newcommand{\bfp}{{\bf p}}%
\newcommand{\bfr}{{\bf r}}%
\newcommand{\bft}{{\bf t}}%
\newcommand{\bfvartheta}{\boldsymbol{\vartheta}}%
\newfont{\tenbfit}{cmbx10}%

\newfont{\tenbbb}{msbm10}%
\newfont{\svnbbb}{msbm8}%

\newcommand{\lj}{\mbox{$[\kern-0.1478125em[$}}
\newcommand{\rj}{\mbox{$]\kern-0.1478125em]$}}
\newcommand{\la}{\mbox{$\langle\kern-0.2325em\langle$}}
\newcommand{\ra}{\mbox{$\rangle\kern-0.2325em\rangle$}}
\newcommand{\Blj}{\mbox{$\Big[\kern-0.275em\Big[$}}
\newcommand{\Brj}{\mbox{$\Big]\kern-0.275em\Big]$}}
\newcommand{\Bla}{\mbox{$\Big\langle\kern-0.425em\Big\langle$}}
\newcommand{\Bra}{\mbox{$\Big\rangle\kern-0.425em\Big\rangle$}}

\newcommand{\zed}{{\bf 0}}

\newcommand{\bfmfrakK}{{\boldsymbol{\mathfrak{K}}}}
\newcommand{\bfmfrakM}{{\boldsymbol{\mathfrak{M}}}}

\begin{document}

\title{Translation and interpretation of Michael Sadowsky's paper ``Theory of elastically bendable inextensible bands with applications to the \textsc{M\"obius} band"\footnote{Citations of this translation should refer also to Sadowsky's original paper, as cited in the Abstract.}
}

\titlerunning{Theory of elastically bendable inextensible bands}        

\author{Denis F.\ Hinz         \and
        Eliot Fried 
}


\institute{Denis F. Hinz \at
              Mathematical Soft Matter Unit\\
              Okinawa Institute of Science and Technology \\
              Okinawa, Japan 904-0495\\
              \email{dfhinz@gmail.com}           
           \and
           Eliot Fried \at
              Mathematical Soft Matter Unit\\
              Okinawa Institute of Science and Technology \\
              Okinawa, Japan 904-0495\\
              \email{eliot.fried@oist.jp}
}

\date{Received: date / Accepted: date}

\maketitle

\begin{abstract}
This article is a translation of Michael Sadowsky's original paper ``Theorie der elastisch biegsamen undehnbaren B\"ander mit Anwendungen auf das \textsc{M\"obius}'sche Band" in 3.\ internationaler Kongress f\"ur technische Mechanik,  Stockholm, 1930.  The translation is augmented by an Appendix containing an interpretation of the last section of Sadowsky's original paper including figures generated from recent numerical simulations.   
\end{abstract}


\section*{Translation of the original paper}

\section{Definition of the term \emph{band} through kinematic properties.}
\label{sec:kinematic}

A body with distinguished midline that, at every point, behaves isotropically with respect to bending is termed a \emph{rope} or a \emph{wire} (or, alternatively, a one-dimensional body). The expression \emph{isotropy with respect to bending} can be illustrated with the following example: If an initially straight and non-twisted section of rope or wire with elastic potential (that is, a conservative system) is bent circularly, then its bending moment depends only on the curvature of the circle and not on which rope or wire elements form the \emph{outer} and \emph{inner} fibers. Consequently, in such a rope or a wire no moment acts to counter eversion of the bent shape.

The circumstances for a \emph{band} stand in contrast to those for a wire or rope. As an example of a band, imagine a \textsc{M\"obius} band made of a thin steel sheet. The surface of the band is a strictly (that is, exactly) developable surface.  The initially planar steel sheet may thus only be bent but not stretched.\footnote{On the questions of existence of an exactly developable \textsc{M\"obius} band, see Sadowsky~\cite{Sadowsky1930}.} A band, no matter how broad or infinitesimally narrow, will respond differently to bending than will a isotropic rope or a wire. In particular, in a band there exists a inadmissible bending axis for relative cross-sectional torsion whereas in a rope or wire all directions have identical properties. In this sense the isotropic response of a rope or wire to bending was mentioned above and in the same sense, a band is anisotropic under bending. A band which has a plane rectangular shape in a stress-free configuration (for example, a \textsc{M\"obius} band) may be bent into circular shape but only in a way that the plane of the circle is everywhere cut perpendicularly by the band. The band thus lies on a cylindrical surface. An inside out rotation of $90^\circ$ by which the band surface would have to take the shape of a part of a circular annulus is impossible because a rectangular surface is not developable on a part of an annulus. Consequently, a band is subject to an eversional moment not present in an isotropic rope or wire.

The distinguishing features of isotropic ropes or wires and bands, approached intuitively here, become clear upon characterizing their midlines and the accompanying vector triads. The three vectors of such a triad will be denoted by $\bft$, $\bfn$, and $\bfb$.

\subsection{Rope or wire}

We fundamentally require that the midline consists of the same particles at all times. It can thus be said that \emph{the midline is material}. Since the tangent vector $\bft$ connects two infinitesimally neighbouring points on the midline, it is firmly connected with the material midline and, thus, we will term it \emph{material} as well. The normal and binormal vectors $\bfn$ and $\bfb$ are unimportant here, because we require the rope or wire to be isotropic with respect to bending. The two latter vectors are thus not associated with the constitution of the rope or wire. Therefore, we may express the fundamental kinematical property of a rope or wire as follows:
\\
\\
\emph{The tangent vector $\bft$ of the accompanying vector triad of a rope or wire is material. }

\subsection{Band}

A band is part of the rectifying surface of its midline. (The rectifying plane of a point $P$ on the midline is the plane through the tangent vector $\bft$ and the binormal vector $\bfb$, the rectifying surface is the envelope of all rectifying planes; cf.\ the work by the author cited in the footnote on the present page). It is simultaneously required that the midline of the band appears as a straight line upon development of the band to a planar figure (namely a non-stressed configuration of the band); in other words, it is required that the midline coincides with a geodesic line on the band. The binormal vector $\bfb$ therefore always lies in the tangent plane of the band, which, however, is \emph{material} in the above mentioned sense---that is firmly connected with the matter of the band. Naturally, the tangent vector $\bft$ remains material as well, since the band is only a special case of an anisotropic rope or wire. In that case, the binormal vector $\bfb$ must be material as well: it lies in a material plane perpendicular to a material direction. Consequently the vectors $\bft$ and $\bfb$ are both material. Moreover, the normal vector 
$$
\bfn = \bfb \times \bft
$$
is material as well, since it is uniquely determined trough the two material vectors $\bft$ and $\bfb$. Thus, the complete vector triad is material. This characteristic kinematical property of a band is encompassed by the following sentence:
\\
\\
\emph{For a band, the accompanying vector triad $\bft$, $\bfn$, and $\bfb$ of the midline is material.}

{\remark{
The statement \emph{the vector triad is material} should not be misconstrued as \emph{the vector triad is made of the same material particles} (which would anyway be senseless for the normal), but should rather be understood as \emph{the vector triad is firmly connected to the matter of the band}.}}

This kinematical property, together with the requirement that the midline is inextensible (that is, the requirement $\delta(ds)=0$) will be viewed as the \emph{definition of a band}.

\section{Conclusions from the kinematical definition of a band: Determination of the virtual torsion of the accompanying triad $\bft$, $\bfn$, and $\bfb$ of the midline compatible with the kinematical properties of the band}

%

\begin{figure}[!t]
\includegraphics[width=0.5\textwidth]{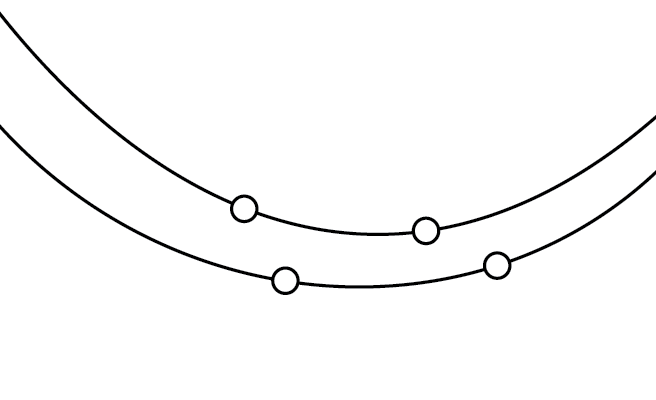}
\put(-108,58){\large$P$}
\put(-56,50){\large$P^\prime$}
\put(-100,16){\large$Q$}
\put(-40,22){\large$Q^\prime$}
 \caption{Adaptation of Figure~1 from the original version of the paper.}
 \label{fig:01}
\end{figure}

\begin{figure}[!t]
\includegraphics[width=0.5\textwidth]{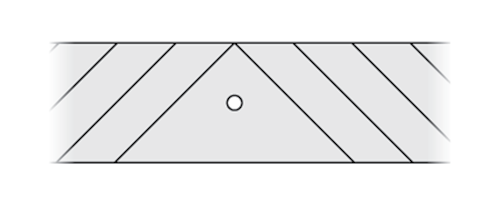}
\put(-86,23){\large$X$}
 \caption{Adaptation of Figure~2 from the original version of the paper.}
 \label{fig:02}
\end{figure}

Let $P$ be an arbitrary point on the midline, $\bfr$ the position vector of that point, $s$ the arclength along the midline, $K$ its curvature (at the point $P$), and $W$ its torsion (at the same point). Let the equation describing the midline be 
\begin{equation}
\bfr = \bfr(s).
\end{equation}

The Frenet--Serret formulas  from the geometry of space curves are then
\begin{equation}\label{eq:frenet01}
\left.
\begin{array}{l}
\displaystyle
\frac{d \bft}{ds} = K \bfn,
\cr\noalign{\vskip8pt}
\displaystyle
\frac{d \bfn}{ds} = - K \bft + W \bfb,
\cr\noalign{\vskip8pt}
\displaystyle
\frac{d \bfb}{ds} = - W \bfn.
\end{array}
\!\!
\right\}
\end{equation}

Let $\delta \bfvartheta$ denote the virtual twist of the accompanying vector triad of the midline at the point $P$. The quantity $\delta \bfvartheta$ should be interpreted as an infinitesimal vector, dependent on the arclength $s$. We now seek to determine the general virtual twist $\delta \bfvartheta$ compatible with the kinematic condition discussed in Section \ref{sec:kinematic}. 

To this end, we consider two neighbouring points $P$ and $P^\prime$ on the midline. The distance between these two points is denoted by $ds$. Since the band is capable only of undergoing bending without stretching, this distance is invariant under the virtual deformation. If the band is subject to an admissible (that is, compatible with the area and band conditions) virtual deformation, then the points $P$ and $P^\prime$ will undergo a displacement taking, say, $P$ to $Q$ and $P^\prime$ to $Q^\prime$. However, the distance $QQ^\prime$ between $Q$ and $Q^\prime$ will remain equal to that $ds$, between $P$ and $P^\prime$ (Figure~\ref{fig:01}). 

Upon application of the virtual deformation, the entire shape of the band changes; specifically, at each point along the midline, the vectors of the triad, as well as the curvature $K$ and the torsion $W$, change. Let
\begin{equation}\label{eq:KW01}
K + \delta K, \quad \text{and} \quad W+\delta W
\end{equation}
denote the curvature and torsion after application of the virtual deformation (both evaluated at the point $Q$, which originates from point $P$ through the virtual deformation). Letting
\begin{equation}
\bft, \quad \bfn, \quad \text{and} \quad \bfb
\end{equation}
denote the accompanying vector triad at the point $P$, we seek to determine the accompanying vector triads at the points $P^\prime$, $Q$, and $Q^\prime$.

To get from $P$ to $P^\prime$, hold the midline fixed and proceed along the curve. For this traversal the Frenet--Serret formulas~\eqref{eq:frenet01} hold and, consequently, we obtain the expressions
\begin{equation}\label{eq:tripod01}
\left.
\begin{array}{l}
\displaystyle
 \bft + K \bfn ds,
\cr\noalign{\vskip8pt}
\displaystyle
 \bfn + (- K \bft + W \bfb) ds,
\cr\noalign{\vskip8pt}
\displaystyle
 \bfb - W \bfn ds,
\end{array}
\!\!
\right\}
\end{equation}
for the three vectors of the triad at the point $P^\prime$.

The point $Q$ is reached after imposing the virtual deformation. For each band element, this deformation may be decomposed into two parts, one involving a virtual translation of the point $P$ and the other involving a virtual rotation around that same point. However, since a virtual translation leaves the vector triad unchanged, only the virtual rotation alters the vectors of the triad. Consequently, to determine the new triad, we need only to consider the virtual rotation $\delta \bfvartheta$. Under an infinitesimal rotation $\delta \bfvartheta$, an arbitrary vector $\bfp$ transforms to the vector $\bfp + \delta \bfvartheta \times \bfp$. For example, the vector $\bft$ becomes the vector $\bft +  \delta \bfvartheta \times \bft$. However, following the kinematic definition of a band in Section~\ref{sec:kinematic}, the vector triad of the midline of a band is {\emph{material}}, that is upon imposing the virtual rotation, the original $\bft$-vector translates into a new $\bft$-vector, the original $\bfb$-vector translates into a new $\bfb$-vector, and so on. We now use this principle and write the following expression for the triad at the point $Q$:
\begin{equation}\label{eq:tripod02}
\left.
\begin{array}{l}
\displaystyle
 \bft + \delta \bfvartheta \times \bft,
\cr\noalign{\vskip8pt}
\displaystyle
 \bfn + \delta \bfvartheta \times \bfn,
\cr\noalign{\vskip8pt}
\displaystyle
 \bfb + \delta \bfvartheta \times \bfb.
\end{array}
\!\!
\right\}
\end{equation}

The point $Q^\prime$ may be reached in two ways. First, consider doing so from $P^\prime$. This involves a transition corresponding to the application of a virtual rotation. In principle, conditions identical to those involved in the transition from $P$ to $Q$ hold; one just needs to be careful to use the vectors corresponding to the triad at $Q^\prime$.  These vectors are given through the formulas~\eqref{eq:tripod01}. Further, the virtual rotation at $P^\prime$ is no longer $\delta \bfvartheta$, but rather is $\delta \bfvartheta + d\delta \bfvartheta$, since the point under consideration is advanced by a distance $ds$ along the midline. On the basis of these considerations, the accompanying vector triad at the point $P^\prime$ becomes
\begin{equation}\label{eq:tripod03}
\left.
\begin{array}{l}
\displaystyle
 \bft + K \bfn ds + (\delta \bfvartheta + d \delta \bfvartheta)\times(\bft + K  \bfn ds),
\cr\noalign{\vskip8pt}
\displaystyle
 \bfn + (- K \bft + W \bfb) ds + (\delta \bfvartheta + d \delta \bfvartheta)\times(\bfn + (- K \bft + W \bfb) ds ) ,
\cr\noalign{\vskip8pt}
\displaystyle
 \bfb - W \bfn ds + (\delta \bfvartheta + d \delta\bfvartheta) \times (\bfb - W \bfn ds).
\end{array}
\!\!
\right\}
\end{equation}

The other way to reach $Q^\prime$ goes through $Q$. This is a displacement along the fixed (virtually displaced) midline. Thus, the Frenet--Serret formulas~\eqref{eq:frenet01} hold for this transition. One only needs to keep in mind that the vector triad at $Q$, as given in~\eqref{eq:tripod02}, corresponds to the starting point and that the curvature and twist of the midline now take the values 
\begin{equation*}
K + \delta K, \quad \text{and} \quad W+\delta W
\end{equation*}
provided in~\eqref{eq:KW01}. The computation gives the vector triad
\begin{equation}\label{eq:tripod04}
\left.
\begin{array}{l}
\displaystyle
 \bft + \delta \bfvartheta \times \bft + (K+\delta K) (\bfn + \delta \bfvartheta \times \bfn) ds,
\cr\noalign{\vskip8pt}
\displaystyle
 \bfn + \delta \bfvartheta \times \bfn+ [- (K+\delta K) \bft + \delta \vartheta \times t + (W+\delta W) (\bfb + \delta \bfvartheta \times \bfb)] ds,
\cr\noalign{\vskip8pt}
\displaystyle
 \bfb + \delta \bfvartheta \times \bfb - (W+
 \delta W)  ( \bfn + \delta \bfvartheta \times \bfn) ds,
\end{array}
\!\!
\right\}
\end{equation}
for the point $Q^\prime$. We now have two representations for the triad in the point $Q^\prime$, namely equations \eqref{eq:tripod03} and \eqref{eq:tripod04}. Consequently, the two sets of equations must describe the same triads of vectors. Through pairwise comparison of the corresponding vectors in~\eqref{eq:tripod03} and~\eqref{eq:tripod04}, one arrives, after neglecting all terms of third order in infinitesimally small quantities, at three equations relating infinitesimally small quantities of second order. These equations read
\begin{equation}\label{eq:relation01}
\left.
\begin{array}{l}
\displaystyle
 \bfn \delta K ds = d \delta \bfvartheta \times \bft,
\cr\noalign{\vskip8pt}
\displaystyle
 -\bft \delta K ds + \bfb \delta W ds =d \delta \bfvartheta \times \bfn ,
\cr\noalign{\vskip8pt}
\displaystyle
 - \bfn \delta W ds = d \delta \bfvartheta \times \bfb .
\end{array}
\!\!
\right\}
\end{equation}
The variable $d\delta \bfvartheta$ appearing in~\eqref{eq:relation01} is the virtual rotation of the vector triad in $P^\prime$ relative to the vector triad in $P$. The system of equations~\eqref{eq:relation01} allows for the determination of the components of $d\delta \bfvartheta$. The relationship 
\begin{equation}\label{eq:sol01}
d \delta \bfvartheta = (\bft \delta W+ \bfb \delta K) ds
\end{equation}
is a general solution to~\eqref{eq:relation01}. {\emph {Together with the constraint $\delta d s = 0$, \eqref{eq:sol01} can be recognized to provide a foundation for the theory of bands.}}

\section{The equations of static equilibrium and the differential equations of a \textsc{M\"obius} band}

The considerations of the two preceding paragraphs hold for bands of arbitrary (that is, not necessarily infinitesimally narrow) width which have a planar shape in the absence of stress and a midline representing a geodesic line on the band surface. From now on,  discussion is restricted to {\emph{infinitesimally narrow}} bands.

Consider a planar cut through the point $P$ on the midline of such a band. The sectional plane shall be the normal plane of the band at the point $P$. The stresses appearing in the sectional plane, which act from the part of the band with larger $s$ to the part of the band with smaller $s$, shall be represented by a force and a moment, both acting at the point $P$. Using the accompanying vector triad at $P$ as a basis for decomposing the force and moment, we write
\begin{equation}\label{eq:force01}
\bfmfrakK = T \bft + N \bfn + B \bfb
\end{equation}
for the force and 
\begin{equation}\label{eq:moment01}
\bfmfrakM = {\mathfrak{T}} \bft + {\mathfrak{N}} \bfn + {\mathfrak{B}} \bfb
\end{equation}
for the moment.

With the requirement that no external forces act on the band, the equilibrium conditions for the band are
\begin{equation}\label{eq:equilForce01}
\frac{ d \bfmfrakK}{ds} = \zed
\end{equation}
and 
\begin{equation}\label{eq:equilForce01}
\frac{ d \bfmfrakM}{ds} +  \bft \times \bfmfrakK = \zed.
\end{equation}
Using the decomposition~\eqref{eq:force01} and~\eqref{eq:moment01}, differentiating in accord with the Frenet--Serret formulas~\eqref{eq:frenet01}, and decomposing the resulting relations into components yields six scalar equations for the six scalar components $T$, $N$, $B$, $\mathfrak{T}$, $\mathfrak{N}$, and $\mathfrak{B}$. The midline curvature $K$ as well as the midline twist $W$ appear in these equations, which read
\begin{equation}\label{eq:relation02}
\left.
\begin{array}{l}
\displaystyle
 \frac{dT}{ds} - KN = 0,
\cr\noalign{\vskip8pt}
\displaystyle
 KT+\frac{dN}{ds} - WB = 0,
\cr\noalign{\vskip8pt}
\displaystyle
 WN+\frac{dB}{ds} = 0,
\cr\noalign{\vskip8pt}
\displaystyle
 \frac{d {\mathfrak T}}{ds} -K {\mathfrak{N}} = 0,
 \cr\noalign{\vskip8pt}
\displaystyle
 K {\mathfrak{T}} + \frac{d {\mathfrak N}}{ds} -W{\mathfrak{B}} -B = 0,
 \cr\noalign{\vskip8pt}
\displaystyle
 W {\mathfrak{N}} + \frac{d {\mathfrak B}}{ds} +N = 0.
\end{array}
\!\!
\right\}
\end{equation}
Notice that \eqref{eq:relation02} do not include a connection between deformation and stress. Consequently, they hold for arbitrary bodies with pronounced midline (rigid or elastic ropes, wires, bands, etc.). To attack our problem, which relates to the equations of the band, from this perspective, we must provide a relationship between deformation and stress. To achieve this, we invoke and follow the principle of virtual displacements.

The theory of virtual work of the internal forces in a one-dimensional continuum with material midline can be found in the literature~\cite{Hamel1926, Hamel1927, Hamel1927a, Hamel1922}. In general, the shear force is always a reaction force, that is, the virtual work of the shear force vanishes identically. Further, for inextensible ropes --- and the band is imagined to be of such a nature --- the  tensile force is a reaction force as well and the corresponding virtual work thus also vanishes identically. The remaining contribution to the internal virtual work $\delta A_i$ is the virtual work of the moment $\bfmfrakM$, which is
\begin{equation}\label{eq:innerWork01}
\delta A_i = - \int \limits _{s_1}^{s_2} \bfmfrakM \cdot d \delta \bfvartheta,
\end{equation}
where $\delta \bfvartheta$ is the virtual rotation of the accompanying vector triad of the midline of the band. However, in~\eqref{eq:sol01} of the preceding paragraph, we  found the most general form for $d\delta \bfvartheta$ that is compatible with the properties of a band. If we now use~\eqref{eq:sol01} together with the component representation~\eqref{eq:moment01} of $\bfmfrakM$, and plug it into~\eqref{eq:innerWork01}, we arrive at an expression,
\begin{equation}\label{eq:innerWork02}
\delta A_i = - \int \limits _{s_1}^{s_2} ({\mathfrak T} \delta W + {\mathfrak B} \delta K) \, ds,
\end{equation}
for the virtual work of the internal forces. {\emph{We thus see that, for a band, the moment $\bfmfrakM$ is also a reaction moment.}}

The relationship~\eqref{eq:innerWork02} holds for every object that may be called \emph{band} in the sense of the definition of bands in Section~\ref{sec:kinematic}. The result~\eqref{eq:innerWork02} is a general relation that holds for arbitrary bands. Now, an arbitrary band is a rather complex object. In the current work, we do not seek to develop a general theory of bands, but rather to explicitly treat only the \textsc{M\"obius} bands. A \textsc{M\"obius} band is characterized energetically through its elastic potential
\begin{equation}
U = A \frac{(K^2+W^2)^2}{K^2},
\end{equation}
where $A$ is a positive constant (see Sadowsky~\cite{Sadowsky1930}, in which a justification for this elastic potential is provided). The elastic energy $E$ of a band segment from $s=s_1$ to $s=s_2$ is then given by
\begin{equation}\label{eq:engery01}
E =  \int \limits _{s_1}^{s_2} U \, ds = A \int \limits _{s_1}^{s_2}  \frac{(K^2+W^2)^2}{K^2} \, ds,
\end{equation}
and the variation of the energy, caused by the virtual deformation of the band, is given as
\begin{equation}\label{eq:varEngery01}
\delta E  = \int \limits _{s_1}^{s_2}  \left( \frac{\partial U}{\partial W} \delta W +\frac{\partial U}{\partial K} \delta K \right) \, ds.
\end{equation}

However, since $\delta A_i = - \delta E$, and since $\delta K$ and $\delta W$ are independent variations, comparison of~\eqref{eq:innerWork02} and~\eqref{eq:varEngery01} leads to the conclusion that
\begin{equation}\label{eq:T01}
{\mathfrak T} = A \frac{4W (K^2+W^2)}{K^2},
\end{equation}
\begin{equation}\label{eq:B01}
{\mathfrak B} = A \frac{2 (K^4-W^4)}{K^3},
\end{equation}
with $A>0$ constant.

The six equations~\eqref{eq:relation02} and the two equations~\eqref{eq:T01} and~\eqref{eq:B01} form a system of eight equations for eight unknowns
\begin{equation*}
T,\, N,\, B,\, {\mathfrak T},\, {\mathfrak N},\, {\mathfrak B},\, K,\, \text{and} \, W.
\end{equation*}
We shall refer to these eight equations as the differential equations of a \textsc{M\"obius} band, since through the integration of these equations one may determine the curvature $K$ and the twist $W$ as a function of the arc length $s$. This then is what may be called a natural description of the midline of the band, and this description determines the unique space curve of the midline of the band up to arbitrary rigid motions. Naturally, these are differential equations that need to be supplemented with appropriate boundary conditions at $s=0$ and $s=l$ (where $l$ is the length of the \textsc{M\"obius} band); these boundary conditions shall account for the correct closure of the strip to a \textsc{M\"obius} band.

Of course, one could eliminate the six coefficients from~\eqref{eq:force01} and~\eqref{eq:moment01} from the eight equations and in the end arrive at two equations for the two geometric variables $K$ and $W$ (this elimination is possible without further difficulties), however the result becomes more complicated as a result, and compromises clarity to such an extent, that this elimination is not conducted here.

The coefficients in~\eqref{eq:force01} and~\eqref{eq:moment01} may be expressed rather clearly through $K$, $W$, and their derivatives with respect to arc length.

The expressions for the coefficients are:
\begin{equation*}
T = AC-U = A\left( C - \frac{(K^2+W^2)^2}{K^2} \right),
\end{equation*}
where $C$ is an arbitrary integration constant,
\begin{equation}\label{eq:dyn02}
N = - \frac{A}{K} \frac{d}{ds} \frac{(K^2+W^2)^2}{K^2},
\end{equation}
\begin{equation}\label{eq:dyn02}
B = \frac{2AW}{K^3}(K^2+W^2)^2 +4A \frac{d}{ds} \left( \frac{1}{K} \frac{d}{ds} \frac{W (K^2+W^2)}{K^2}\right),
\end{equation}
\begin{equation}\label{eq:dyn03}
\mathfrak N = \frac{4A}{K}\frac{d}{ds} \frac{W (K^2+W^2)}{K^2},
\end{equation}
\begin{center}($\mathfrak T$ and $\mathfrak B$, cf.~\eqref{eq:T01} and~\eqref{eq:B01}).\end{center}
These six equations further need to satisfy the conditions
\begin{equation}\label{eq:dyn04}
KT + \frac{dN}{ds} - WB = 0
\end{equation}
and
\begin{equation}\label{eq:dyn05}
WN + \frac{dB}{ds} = 0.
\end{equation}
Plugging~\eqref{eq:dyn03} into~\eqref{eq:dyn04} and~\eqref{eq:dyn05} would give the two differential equations for the curvature $K$ and the twist $W$.

\section{A peculiar implication for the shape of a \textsc{M\"obius} band}
\label{sec:peculiar}

Due to a lack of space, this last section can only be presented in the shortest possible form. Following the considerations with the help of a paper model and scissors is therefore recommended, as these considerations might otherwise remain incomprehensible.

The midline of the band possesses a singular point $X$ that can be found in the following way: a \textsc{M\"obius} band possesses a symmetry axis defined such that it is congruent with itself after a rotation of $180^\circ$ about that axis. The symmetry axis intersects the band at two points such that the axis coincides with the binormal $\bfb$ at one of these points. In Figure~\ref{fig:02}, this is the point $X$.

Let $\varphi$ denote the angle between the rectilinear generators of the band through the point $P$ and $\bfb$; then
\begin{equation}\label{eq:phi01}
\tan \varphi = \frac{W}{K}.
\end{equation}

For the point $X$,
\begin{equation}\label{eq:Beq0}
{\mathfrak B} =0
\end{equation}
due to symmetry.

To arrive at an additional conclusion, consider an {\emph{experiment with a band model}}: this shows that
\begin{equation}\label{eq:exp01}
\lim \limits_{P \rightarrow X} \varphi \neq 0,
\end{equation}
\begin{equation}\label{eq:exp02}
\lim \limits_{P \rightarrow X} {\mathfrak T} \neq 0.
\end{equation}
The experiment corresponding to~\eqref{eq:exp01} consists of observing a band; for~\eqref{eq:exp02}, one needs to cut the boundaries of a band at $X$ to infer the moment $\mathfrak T$ from the twist of the band in the weakened cross section. By~\eqref{eq:Beq0},~\eqref{eq:exp01}, \eqref{eq:exp02}, and the governing equations~\eqref{eq:T01} and~\eqref{eq:B01} of the band, it transpires that 
\begin{equation}\label{eq:lim01}
\lim \limits_{P \rightarrow X} K \neq 0, \quad \lim \limits_{P \rightarrow X} W \neq 0, \quad \text{and} \quad \lim \limits_{P \rightarrow X} \varphi = 45^\circ.
\end{equation}
The last limit in~\eqref{eq:lim01} points to the following peculiar fact:\\

{\emph{A \textsc{M\"obius} band consists of a planar, right triangle. The curved, analytic portion of the band connects continuously with the two legs of the right triangle; it connects with continuous tangential plane, but with discontinuous curvature.}}


\appendix

\section{Interpretation and explanation of Section~\ref{sec:peculiar} of Sadowsky's paper}

Unfortunately, Section~\ref{sec:peculiar} of Sadowsky's paper is too brief to contain illustrative explanations as to why this ``peculiar implication" is so peculiar and important. In this appendix, we attempt to provide an additional illustrative explanation of Sadosky's early observation. To this end, we discuss the last section in detail and use recent simulation results by Kleiman et al.~\cite{Kleiman2014}.

The curvature of the midline $K$ represents the curvature along the tangent direction of the midline, which cuts through the rectilinear generators of the band surface at an angle. This angle is $90^\circ \pm \varphi$, where $\varphi$ is the angle used by Sadowsky in~\eqref{eq:phi01}. This geometric consideration (or a similar geometric argument) leads to~\eqref{eq:phi01}. 

\begin{figure*}[t]
\hspace{0.9cm}
\includegraphics[width=.9\textwidth] {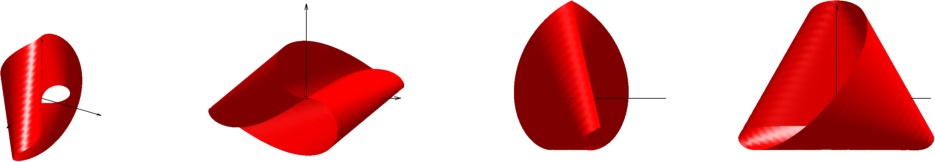}
\put(-325,20){$ a = \pi$}
\put(-34,30){$\color[rgb]{0.7, 0.75, 0.71}\circ$}
\put(-40,30){$X$}
\put(-120,30){$\color[rgb]{0.7, 0.75, 0.71}\circ$}
\put(-126,30){$X$}
\put(-206.5,19){$\color[rgb]{0.7, 0.75, 0.71}\circ$}
\put(-212.5,19){$X$}
\put(-80,55){symmetry axis}
\put(-50,50){\large \rotatebox{-40}{$\longrightarrow$} }
\put(-307,8){$x$}
\put(-275,10){$y$}
\put(-289,43){$z$}
\put(-175,15){$x$}
\put(-202,48){$y$}
\put(-90,15){$y$}
\put(-115,48){$z$}
\put(-3,15){$x$}
\put(-30,48){$z$}

\hspace{0.9cm}
\includegraphics[width=.9\textwidth] {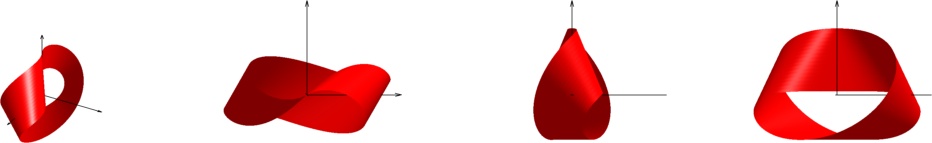}
\put(-325,20){$ a = 2\pi$}
\put(-33,25){$\color[rgb]{0.7, 0.75, 0.71}\circ$}
\put(-40,25){$X$}
\put(-120,25){$\color[rgb]{0.7, 0.75, 0.71}\circ$}
\put(-126,25){$X$}
\put(-206,14){$\color[rgb]{0.7, 0.75, 0.71}\circ$}
\put(-212,14){$X$}
\put(-307,3){$x$}
\put(-275,5){$y$}
\put(-290,34){$z$}
\put(-175,10){$x$}
\put(-202,43){$y$}
\put(-90,10){$y$}
\put(-115,43){$z$}
\put(-3,10){$x$}
\put(-30,43){$z$}

\hspace{0.9cm}
\includegraphics[width=.9\textwidth] {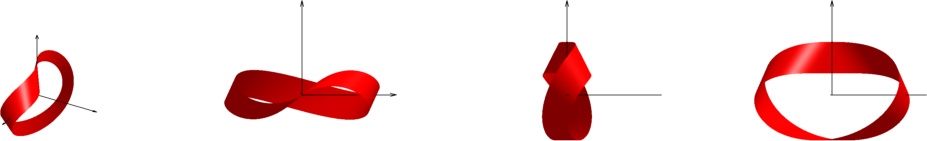}
\put(-325,20){$ a = 4\pi$}
\put(-33,25){$\color[rgb]{0.7, 0.75, 0.71}\circ$}
\put(-40,25){$X$}
\put(-120,25){$\color[rgb]{0.7, 0.75, 0.71}\circ$}
\put(-126,25){$X$}
\put(-207,13){$\color[rgb]{0.7, 0.75, 0.71}\circ$}
\put(-213,13){$X$}
\put(-307,3){$x$}
\put(-275,5){$y$}
\put(-290,34){$z$}
\put(-175,10){$x$}
\put(-202,43){$y$}
\put(-90,10){$y$}
\put(-115,43){$z$}
\put(-3,10){$x$}
\put(-30,43){$z$}

\hspace{0.9cm}
\includegraphics[width=.9\textwidth] {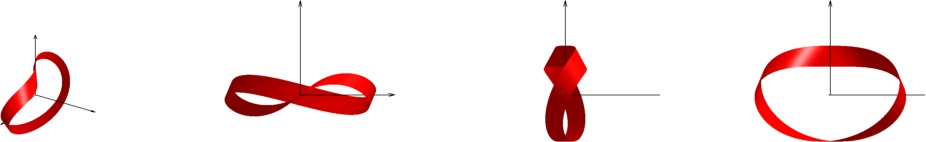}
\put(-325,20){$ a = 6\pi$}
\put(-33.5,26){$\color[rgb]{0.7, 0.75, 0.71}\circ$}
\put(-40,26){$X$}
\put(-120,26){$\color[rgb]{0.7, 0.75, 0.71}\circ$}
\put(-126,26){$X$}
\put(-207.5,13.5){$\color[rgb]{0.7, 0.75, 0.71}\circ$}
\put(-213.5,13.5){$X$}
\put(-307,3){$x$}
\put(-275,5){$y$}
\put(-290,34){$z$}
\put(-175,10){$x$}
\put(-202,43){$y$}
\put(-90,10){$y$}
\put(-115,43){$z$}
\put(-3,10){$x$}
\put(-30,43){$z$}

\hspace{0.9cm}
\includegraphics[width=.9\textwidth] {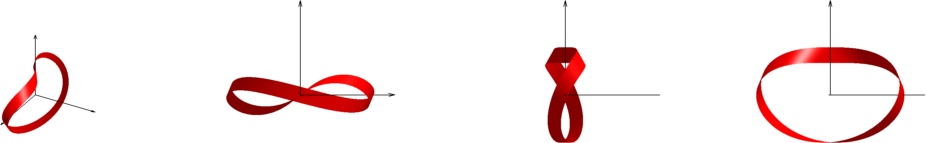}
\put(-325,20){$ a = 8\pi$}
\put(-33.5,27){$\color[rgb]{0.7, 0.75, 0.71}\circ$}
\put(-40,27){$X$}
\put(-120,27){$\color[rgb]{0.7, 0.75, 0.71}\circ$}
\put(-126,27){$X$}
\put(-207.5,13.5){$\color[rgb]{0.7, 0.75, 0.71}\circ$}
\put(-213.5,13.5){$X$}
\put(-307,3){$x$}
\put(-275,5){$y$}
\put(-290,34){$z$}
\put(-175,10){$x$}
\put(-202,43){$y$}
\put(-90,10){$y$}
\put(-115,43){$z$}
\put(-3,10){$x$}
\put(-30,43){$z$}

\caption{Location of the singular point $X$ and the (rotational) symmetry axis: Equilibrium shapes of approximately developable M\"obius bands for different aspect ratios $a$ obtained from simulations with a lattice model by Kleiman et al.$^2$ The band is rotated into its main axes and the approximate location of the point $X$ is indicated. With the current choice of coordinate system, the (rotational) symmetry axis coincides with the $z$-axis. }
\label{fig:Moebius_shape}
\end{figure*}

Further,  the (rotational) symmetry axis of a M\"obius band can indeed easily be found with a paper model. In the simulation results, the symmetry axis coincides with the $z$-axis of the coordinate system, as shown in Figure~\ref{fig:Moebius_shape}. Due to symmetry, the bending moment ${\mathfrak B}$ in~\eqref{eq:B01} vanishes at $X$ (see~\eqref{eq:Beq0}). Using the condition~\eqref{eq:Beq0} in~\eqref{eq:B01} yields 
\begin{equation}
K^4-W^4=0
\end{equation}
and, thus, 
\begin{equation}\label{eq:KpmW}
K =\pm W
\end{equation}
at $X$, which includes the possibility that $K=0$ and $W=0$ at $X$.

\begin{figure*}[t]
\begin{center}
\includegraphics[width=.8\linewidth] {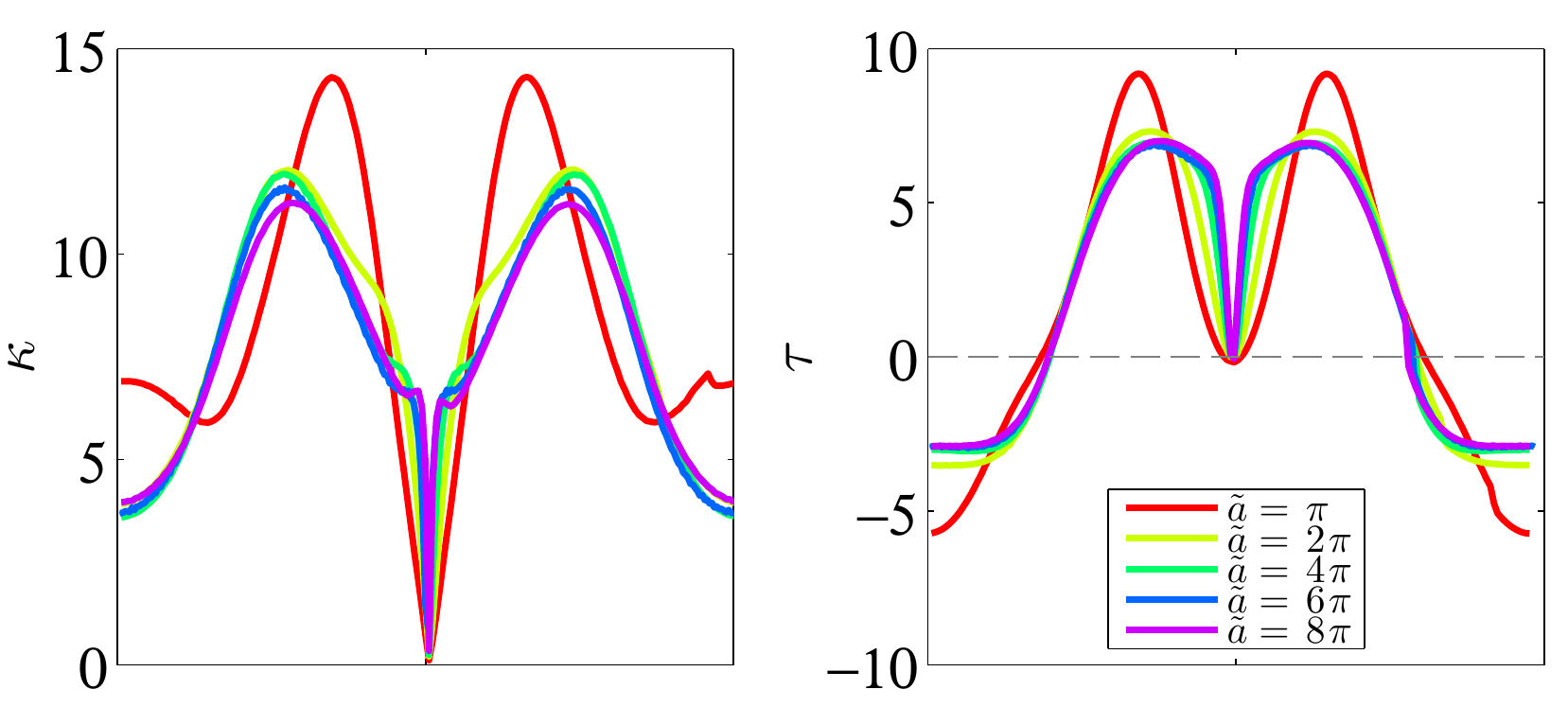}
\put(-81,8){\includegraphics[width=.14\linewidth] {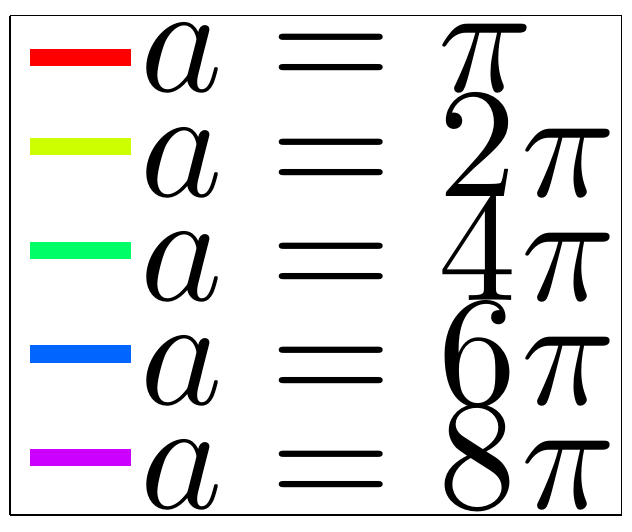}}
\put(-59.5,58.5){$\circ$}
\put(-59.5,50){$X$}
\put(-198.5,5){$\circ$}
\put(-205,9){$X$}
\put(-60,-13){\fontsize{0.4cm}{1em}\selectfont$s$}
\put(-200,-13){\fontsize{0.4cm}{1em}\selectfont$s$}
\put(-113,-3){\large$0$}
\put(-60,-3){\large$\pi$}
\put(-10,-3){\large$2\pi$}
\put(-253,-3){\large$0$}
\put(-200,-3){\large$\pi$}
\put(-150,-3){\large$2\pi$}
\put(-136,55){\fcolorbox{white}{white}{\rotatebox{90}{\fontsize{0.4cm}{1em}\selectfont$W$}}}
\put(-276,55){\fcolorbox{white}{white}{\rotatebox{90}{\fontsize{0.4cm}{1em}\selectfont$K$}}}

\end{center}
\caption{The curvature and twist (in arbitrary units) of developable M\"obius bands show the singular point $X$ with $K=0$ and $W=0$: Curvature and twist of the centerline of the developable M\"obius bands for different aspect ratios $a$. The singular nature of this point becomes more prominent for increasing aspect ratio $a$. Results are obtained from simulations with a lattice model by Kleiman et al.$^2$ where complete results, validation, and discussion can be found. }
\label{fig:curvatureTwist}
\end{figure*}

Now, suppose that $K\neq 0$ at $X$. In view of~\eqref{eq:KpmW},
\begin{equation}
\frac{W}{K}=\pm 1
\end{equation}
and, thus,
\begin{equation}\label{eq:tanphi}
\tan \varphi=\pm 1.
\end{equation}
The angle $\varphi$ would thus be
\begin{equation}\label{eq:45deg}
\varphi = \pm \frac{1}{4}\pi = \pm 45^\circ.
\end{equation}
Equation~\eqref{eq:45deg} implies that there must be two perpendicular rectilinear generators of the bent surface at $X$, which means that the bent surface must be flat at $X$. However, a flat surface at $X$ contradicts the supposition that $K\neq0$ at $X$. Consequently, 
\begin{equation}\label{eq:KeqWeq0}
K = 0 \quad {\text{and}} \quad W = 0
\end{equation}
at $X$. Note that the rectilinear generator at $X$ coincides with $\bfb$ at $X$ by symmetry. The result~\eqref{eq:KeqWeq0} is confirmed by plots generated by our numerical simulations shown in Figure~\ref{fig:curvatureTwist}. 

Now, consider the point $P$ on the midline close to $X$. The curvature $K$ cannot be zero at this point. This becomes clear by the following observation. If $K$ was zero at $P$, the surface would have zero curvature along the rectilinear generator through $P$ as well as along the midline tangent direction through $P$. Since these two tangent directions are not equal, the surface would have to be flat at $P$, and not bent, as is presumed. Consequently, 
\begin{equation}\label{eq:kneq0}
K\neq 0
\end{equation}
at $P$. With~\eqref{eq:kneq0} holding at $P$, $W$ may or may not vanish. If $W$ vanishes at $P$, then the rectilinear generator through $P$ would have to coincide with $\bfb$ at $P$, since $\varphi=0$ by~\eqref{eq:phi01}.  That means that, for $W=0$, the rectilinear generator at P is perpendicular to the tangent direction of the midline at $P$, in which case $K$ represents the single non-zero principal curvature of the bent surface at $P$. In addition, for $W=0$, in view of~\eqref{eq:T01}, there would be no twisting moment around the axis tangent to the midline at $P$. However, this condition of $K\neq 0$ and $W=0$ cannot persist for all $P$ on the midline, since such a configuration would not be consistent with an isometric mapping to a bent (and twisted) M\"obius band. In other words, for $K$ to vanish, $W$ must vanish as well (corresponding to the singular point $X$), whereas $W$ may vanish without $K$ vanishing, but not for all points. This is confirmed by our simulation results (Figure~\ref{fig:curvatureTwist}), which consistently show three zeros for $W$ and non-zero $W$ on the rest of the band. 

In view of the preceding discussion, we conclude that what Sadowsky was trying to communicate with~\eqref{eq:exp01}, \eqref{eq:exp02}, and~\eqref{eq:lim01} is that any point $P$ on the midline in a neighborhood of $X$ has $K \neq 0$ and must also have $W \neq 0$, which, according to~\eqref{eq:T01}, ensures that there is a non-zero twisting moment about the tangent axis through $P$ and, according to~\eqref{eq:phi01}, $\varphi \neq 0$ at $P$. Moreover, the condition $W \neq 0$ persists up to the point $X$, where both $K$ and $W$ are zero. In other words, since $K=0$ and $W=0$ at $X$ and $K\neq 0$ and $W\neq 0$ at $P$, the generators connect to the flat point $X$ at $\pm 45^\circ$. This corresponds to the ``planar, right triangle" to which Sadowsky refers. In fact, such planar, right triangles are well observed in the simulation results shown in Figure~\ref{fig:Moebius_shape}, especially for low aspect ratios. For high aspect ratios, the singular nature of the curvature and twist at $X$ becomes more prominent, as shown in Figure~\ref{fig:curvatureTwist}. Considering that Sadowsky did not have access to simulation results, his early observation is quite remarkable.

\begin{acknowledgements}

The authors thank Roger L.~Fosdick, Russel E.~Todres, and David M.~Kleiman for enlightening discussions during the translation of Sadowsky's paper, discussions which led to the inclusion of the Appendix for clarification of Sadowsky's original considerations. 

\end{acknowledgements}

%
%
%




\end{document}